\documentclass[12pt]{article}
\usepackage{amsmath,amsfonts,amstext,amsthm,epsfig}
\begin{document}
\font\tenrm=cmr10

\newtheorem{Theorem}{Theorem}
\newtheorem{Lemma}{Lemma}
\newtheorem{Corol}{Corollary}

\newtheorem{Prop}{Proposition}
\newtheorem{Rem}{Remark}
\makeatletter
    \addtocounter{section}{0}
    \renewcommand{\theequation}{\thesection.\arabic{equation}}
    \@addtoreset{equation}{section}
\makeatother

\title{\bf Backward uniqueness for parabolic operators with non--Lipschitz
coefficients}

\author{ Daniele Del Santo\\ Martino Prizzi\\ {\small Dipartimento di
Scienze Matematiche,
Universit\`a di Trieste}\\[-0.1 cm]
 {\small Via A.~Valerio 12/1, 34127
 Trieste, Italy} }
\date{\today}
\maketitle

\begin{abstract} We investigate the relation between the backward
uniqueness and the regularity of
the coefficients for a parabolic operator. A necessary and sufficient
condition for uniqueness is
given in terms of the modulus of continuity of the coefficients.
\vskip 0.3 cm\noindent {\bf Keywords}: backward uniqueness, parabolic
operators, modulus of
continuity, Osgood condition\end{abstract}
\section {Introduction} We consider the following backward parabolic operator
\begin{equation} L=\partial_t+ \sum_{i, j = 1}^n \partial_{x_j} (a_{jk} (t,x)
 \partial_{x_k} )+\sum_{ j = 1}^n b_j(t,x)\partial_{x_j}+c(t,x).\label{1}
\end{equation} All the coefficients are supposed to be defined in $[0,T]
\times {\mathbb R}^n_x$,
measurable  and bounded; the coefficients $b_j$ and $c$ are complex valued;
$( a_{jk}(t,x))_{jk}$ is a real symmetric matrix for
all $(t,x)\in [0,T]\times
{\mathbb R}^n_x$ and there exists $\lambda_0\in (0,1]$ such that
$$
\sum_{j, k = 1}^n  a_{jk} (t,x)
\xi_j \xi_k\geq \lambda_0|\xi|^2
$$ for all $(t,x)\in [0,T]\times  {\mathbb R}^n_x$ and $\xi\in {\mathbb
R}^n_\xi$.

Given a functional space ${\cal H}$ (in which it makes sense to look for
the solutions of the
equation $Lu=0$) we say that the operator $L$ has the ${\cal
H}$--uniqueness property if, whenever $u\in
{\cal H}$, $Lu=0$ in $[0,T]\times  {\mathbb R}^n_x$ and
$u(0,x)=0$ in ${\mathbb R}^n_x$, then $u=0$ in $[0,T]\times  {\mathbb R}^n_x$.

The problem we are interested in is the following:  find the minimal
regularity on the coefficients
$a_{jk}$ ensuring the ${\cal H}$--uniqueness property to $L$.

We remark that even in the simplest case (i. e. $(a_{jk})_{jk}={\rm Id}$) the
answer may depend on
${\cal H}$ and in particular on the rate of growth of $u$ with respect to
the $x$ variables, as the
classical example of Tychonoff  \cite{Ty} shows.

Considering $ {\cal H}_1= H^1([0,T], L^2({\mathbb R}^n_x))\cap L^2([0,T],
H^2({\mathbb R}^n_x))$, $ {\cal
H}_1$--uniqueness for $L$ has been proved under the hypothesis of
Lipschitz--continuity of the
coefficients $a_{jk}$ by Lions and Malgrange \cite {LM} (see for related or
more general results
\cite {Mi}, \cite{AN}, \cite{BT}, \cite{Gh}). On the other hand the well
known example of Miller
\cite{Mil} (where an operator having coefficients which are
H\"older--continuous of order $1/6$
with respect to $t$ and $C^\infty$ with respect to $x$ does not have the
uniqueness property) shows
that a certain amount of regularity on the $a_{jk}$'s is necessary for the
$ {\cal
H}_1$--uniqueness.

The first part of the present work is devoted to prove the $ {\cal
H}_1$--uniqueness property for
the operator (\ref{1}) when the coefficients $a_{jk}$ are $C^2$ in the $x$
variables and
non--Lipschitz--continuous in $t$. The regularity in $t$ will be given in
terms of a modulus of
continuity $\mu$ satisfying the so called Osgood condition
$$
\int_0^1{1\over \mu(s)}\, ds=+\infty.
$$

This uniqueness result is a consequence of a Carleman estimate in which the
weight function depends
on the modulus of continuity;  such kind of weight functions in Carleman
estimates have been
introduced by Tarama \cite {Ta} in the case of second order elliptic
operators. In obtaining our
Carleman estimate the integrations by parts, which cannot be used since the
coefficients are not
Lipschitz--continuous,  are replaced by a microlocal approximation
procedure similar to the one
exploited by Colombini and Lerner \cite{CL} to prove some energy estimates
for hyperbolic operators
with log--Lipschitz coefficients (see also \cite{CDGS} and \cite{CDS}).

It is  interesting to remark that the Osgood condition is also necessary
for the  $ {\cal
H}_1$--uniqueness property, at least when only the regularity in $t$ of the
coefficients $a_{jk}$
is concerned. Precisely in the second part of this paper we prove that if a
modulus of continuity
does not satisfy the Osgood condition then it is possible to construct a
backward parabolic operator
of type (\ref{1}) such that the coefficients $a_{jk}$ depend only on $t$,
the regularity of the
$a_{jk}$'s is ruled by the modulus of continuity and the operator has not
the  $ {\cal
H}_1$--uniqueness property. The construction of this class of examples is
modelled on a well known
non--uniqueness result for elliptic operators due to Pli\'s \cite{Pl}.

The plan of the paper is the following: in Section 2 we give the precise
statement of the
uniqueness theorem and we present the non--uniqueness examples; a remark is
devoted to compare
these results with similar ones known for elliptic and hyperbolic
operators. Section 3 contains
the proof of the uniqueness results. In Section 4 we sketch the
construction of the counter examples.
\par  We denote by $\langle\cdot,\cdot\rangle_{L^2}$ the scalar
product in $L^2({\mathbb R}^n_x)$ and by $\|\cdot\|_{L^2}$ the
corresponding norm. We denote by $\|\cdot\|_{\cal B}$ the norm of
any other  Banach space ${\cal B}$. Finally we denote by $\nabla$
the gradient with respect to the $x$ variables.

\section {Results and remarks}

Let $\mu$ be a modulus of continuity, i. e. let $\mu\,:\, [0,1]\to [0,1]$
be continuous, concave,
strictly increasing, with $\mu(0)=0$. Let $I \subseteq {\mathbb R}$ and
let
$\varphi \,:\, I\to {\cal B}$, where ${\cal B}$ is a Banach space. We say
that $\varphi\in
C^\mu(I,{\cal B})$ if $\varphi\in L^\infty(I,{\cal B})$ and
$$
\sup_{0<|t-s|<1\atop t,s\in I}{\|\varphi(t)-\varphi (s)\|_{\cal
B}\over\mu(|t-s|)}<+\infty.
$$

\begin{Rem} \label{rem1}
The concavity of $\mu$ implies that $\mu(s)\geq s\mu(1) $ for
all $s\in [0,1]$; the
same reason makes the function $s\mapsto {\mu(s)\over s}$ decreasing on
$\;]0, 1]$. Consequently
there exists $\lim _{s\to 0^+} {\mu(s)\over s}$. If $\sup_{s\in \, ]0,1]}
{\mu(s)\over s}<+\infty$
then  there exists $C>0$ such that $\mu(s)\leq Cs$ for all $s\in [0,1]$ and
hence $C^\mu={\rm
Lip}$. As a consequence,  if
$C^\mu\not ={\rm Lip}$, in particular if $\int^1_01/\mu(s)\, ds<+\infty$,
then $\lim _{s\to 0^+}
{\mu(s)\over s}=+\infty$. Finally the function $\sigma\mapsto
\mu(1/\sigma)/(1/\sigma)$ is
increasing on $[1,+\infty[$;  consequently  the function $\sigma\mapsto
\sigma^2\mu(1/\sigma)$ is
increasing on $[1,+\infty[$ and the function
$\sigma\mapsto 1/(\sigma^2\mu(1/\sigma))$ is decreasing on the same interval.
\end{Rem}

We can now state our main uniqueness result.
\begin{Theorem} Let $\mu$ be a modulus of continuity and suppose
\begin{equation}
\int_0^1{1\over \mu(s)}\, ds=+\infty.
\label{3}
\end{equation} Suppose, for all $j,k=1\dots,n$, $a_{jk}\in
C^\mu([0,T],C^2_b({\mathbb R}^n_x))$ where
$C^2_b({\mathbb R}^n_x)$ is the space of twice differentiable functions
which are
bounded with bounded
derivatives.

Then the operator $L$ defined in (\ref{1}) has the ${\cal H}_1$--uniqueness
property.
\label{t1}
\end{Theorem}

Let us denote by ${\cal H}_2$ the space of functions $w$ defined in
$[0,T]\times{\mathbb R}^n_x$ such that
$w$ is continuous and differentiable with respect to $t$ with
continuous derivative and
twice differentiable with respect to $x$ with continuous derivatives and
there exists
$C>0$ such that
$$ |w(t,x)|,\; |\partial_t w(t,x)|,\; |\partial_{x_j}
w(t,x)|,\;|\partial_{x_j}\partial_{x_k}
w(t,x)|\leq Ce^{C|x|}
$$ for all $j,k=1\dots,n$ and for all $(t,x)\in [0,T]\times{\mathbb
R}^n_x$. The
following result holds.

\begin{Theorem} In the hypotheses of Theorem \ref{t1} the operator $L$ has
 the ${\cal H}_2$--uniqueness property.
\label{t2}
\end{Theorem}

The condition (\ref{3}) on $\mu$ is  known as  ``Osgood condition"
(see e.g.
\cite[p. 160]{Fl}. Our next result shows that this condition is necessary
to have the uniqueness
property.

\begin{Theorem} Let $\mu$ be a modulus of continuity and suppose
\begin{equation}
\int_0^1{1\over \mu(s)}\, ds<+\infty.
\label{4}
\end{equation}

Then there exists $l\in C^\mu({\mathbb R}_t)$ with $1/2\leq l(t)\leq 3/2$
for all
$t\in {\mathbb R_t}$ and
there exists $u$, $b_1$, $b_2$, $c\in C^\infty_b({\mathbb R}_t
\times{\mathbb R}^2_x)$ with ${\rm supp}\;
u=\{t\geq 0\}$ such that
\begin{equation}
\partial_t u+ \partial^2_{x_1}u+ l
\partial^2_{x_2}u+b_1\partial_{x_1}u+b_2\partial_{x_2}u+cu=
0\quad {\it in}\ {\mathbb R}_t\times{\mathbb R}^2_x.
\label {5}
\end{equation}
\label{t3}
\end{Theorem}

\begin{Rem} Considering a function $\theta\in C^\infty({\mathbb R}^n_x)$
such that
$\theta(x)=e^{-C|x|}$ for $|x|\geq 1$ and taking $v(t,x)=
\theta(x)u(t,x)$ where $u(t,x)$ is the function constructed in Theorem
\ref{t3}, we immediately
obtain a counter example to the
${\cal H}_1$--uniqueness result.
\end{Rem}

\begin{Rem} It may be interesting to compare the uniqueness and
non--uniquen\-ess results presented
here with similar ones known for different classes of operators. The case
of second order elliptic
operators with real principal part has been considered by Tarama
\cite{Ta}. The uniqueness in the
Cauchy problem is obtained for such kind of operators when the
coefficients of the principal part
are $C^\mu$ with respect to all the variables and $\mu $ satisfies the
condition (\ref{3}). A
precise analysis of the non--uniqueness example of Pli\'s \cite{Pl} shows
that (\ref{3}) is
necessary (see \cite{DS}).

An example of non--uniqueness for hyperbolic operators having the
coefficients of the principal
part in $C^\mu$ with $\mu$ satisfying the condition (\ref{4}) is given in
\cite {CDS2} (see also
\cite{CJS}). It is an open problem, whether (\ref{3}) is sufficient to have the
uniqueness in the Cauchy
problem for second order hyperbolic operators.
\end{Rem}

\section{Proofs of Theorems 1 and  2}

In this paragraph we prove Theorem \ref{t1} and Theorem \ref{t2}. Theorem
\ref{t1} will follow in
standard way from a Carleman estimate. In order to state the latter, we
need first to introduce the weight function. We
define
$$
\phi(t)=\int_{ 1\over t}^1{1\over \mu(s)}\; ds.
$$ The function $\phi$ is a strictly increasing $C^1$ function. From
(\ref{3}) we have
$\phi([1,+\infty[)=[0,+\infty[$ and  $\phi'(t)= 1/(t^2\mu(1/t))>0$ for all
$t\in [1,+\infty[$. We
set
$$
\Phi(\tau)=\int_0^\tau\phi^{-1}(s)\;ds.
$$ We obtain $\Phi'(\tau)=\phi^{-1}(\tau)$ and consequently $\lim_{\tau\to
+\infty}\Phi'(\tau)=+\infty$. Moreover
\begin{equation}
\Phi''(\tau)=(\Phi'(\tau))^2\mu({1\over \Phi'(\tau)})
\label{6}
\end{equation} for all $\tau\in [0,+\infty[$ and, as the function
$\sigma\mapsto
\sigma\mu(1/\sigma)$ is increasing on $[1,+\infty[$ (see Remark \ref{rem1}), we
deduce that
\begin{equation}
\lim_{\tau\to +\infty}\Phi''(\tau)=\lim_{\tau\to
+\infty}(\Phi'(\tau))^2\mu({1\over
\Phi'(\tau)})=+\infty.
\label{6bis}
\end{equation} Now we can state the Carleman estimate.
\begin{Prop} There exist $\gamma_0$, $C>0$ such that
\begin{equation}
\begin{array}{ll}
\displaystyle{\int_0^{T\over 2} e^{{2\over
\gamma}\Phi(\gamma(T-t))}\|\partial_tu+\sum_{jk}
\partial_{x_j}(a_{jk}\partial_{x_k} u) \|^2_{L^2}\;
dt}\\[0.3 cm]
\displaystyle{\qquad\geq C\gamma^{1\over 2}\int_0^{T\over 2} e^{{2\over
\gamma}\Phi(\gamma(T-t))}
(\|\nabla u\|^2_{L^2}+\gamma^{1\over
2}\|u\|^2_{L^2})\; }dt
\end{array}
\label{7}
\end{equation} for all $\gamma>\gamma_0$ and for all $u\in
C^\infty_0({\mathbb R}_t\times{\mathbb R}^n_x, {\mathbb C})$
such that $\,{\rm supp}\, u\subseteq [0, T/2]\times {\mathbb R}^n_x$.
\label{p1}
\end{Prop}

The proof of the Proposition \ref{t1} is rather long and we divide it in
several steps.
\vskip 0.2 cm
\noindent {\it a) the Littlewood--Paley decomposition}
\vskip 0.2 cm
\noindent We set $v(t,x)=e^{{1\over \gamma}\Phi(\gamma(T-t))} u(t,x)$.  The
inequality (\ref{7})
becomes
\begin{equation}
\begin{array}{ll}
\displaystyle{\int_0^{T\over 2} \|\partial_tv+\sum_{jk}
\partial_{x_j}(a_{jk}\partial_{x_k} v) + \Phi'(\gamma(T-t))
v\|^2_{L^2}\;
dt}\\[0.3 cm]
\displaystyle{\qquad\geq C\gamma^{1\over 2}\int_0^{T\over 2}
(\|\nabla v\|^2_{L^2}+\gamma^{1\over
2}\|v\|^2_{L^2})\; }dt.
\end{array}
\label{8}
\end{equation} We use now the Littlewood--Paley decomposition technique.
We recall some basic
facts on it, referring to \cite{Bo} and \cite{CL} for further details. Let
$\varphi_0\in
C^\infty_0({\mathbb R}^n_\xi)$, $0\leq \varphi(\xi)\leq 1$ for all $\xi\in
{\mathbb R}^n_\xi$,
$\varphi_0(\xi)=1$ for all $\xi$ such that $|\xi|\leq 1$,
$\varphi_0(\xi)=0$ for all $\xi$ such
that $|\xi|\geq 2$ and $\varphi_0$ radially decreasing. For all $\nu\in
{\mathbb N}\setminus\{0\}$ we
define
$$
\varphi_\nu(\xi)=\varphi_0({\xi\over 2^\nu})-\varphi_0({\xi\over 2^{\nu-1}}).
$$ For $u\in L^2({\mathbb R}^n_x,{\mathbb C})$ we set
\begin{equation} u_\nu (x)= \varphi_\nu(D)u(x)=
{1\over(2\pi)^n}\int_{{\mathbb R}^n_\xi}
e^{{\rm i}x\xi} \varphi_\nu(\xi)\hat
u(\xi)\; d\xi,
\label{9}
\end{equation}
where $\hat u$ is the Fourier-Plancherel transform of $u$. We remark that
(\ref{9}) makes sense also for $u\in {\cal S}'({\mathbb R}^n_x,{\mathbb
C})$ if the
last integral is interpreted as the inverse Fourier transform of
$\varphi(\xi)\hat u(\xi)$. We have
that there exists $K>0$ such that
\begin{equation} {1\over K}\sum_\nu\|u_\nu\|^2_{L^2}\leq \|u\|^2_{L^2}\leq
K\sum_\nu\|u_\nu\|^2_{L^2}
\label{10}
\end{equation} for all $u\in L^2({\mathbb R}^n_x,{\mathbb C})$. Consequently
\begin{equation}
\begin{array}{l}
\displaystyle{\int_0^{T\over 2} \|\partial_tv+\sum_{jk}
\partial_{x_j}(a_{jk}\partial_{x_k} v) + \Phi'(\gamma(T-t))
v\|^2_{L^2}\;
dt}\\[0.3 cm]
\displaystyle{\qquad\geq {1\over K}\int_0^{T\over 2}
\sum_\nu\|\varphi_\nu(D)(\partial_tv+\sum_{jk}
\partial_{x_j}(a_{jk}\partial_{x_k} v) + \Phi'(\gamma(T-t)) v)\|_{L^2}^2\;
dt}\\[0.3 cm]
\displaystyle{\qquad\geq {1\over K}\int_0^{T\over 2}
\sum_\nu\|\partial_tv_\nu+\sum_{jk}
\partial_{x_j}(a_{jk}\partial_{x_k} v_\nu) + \Phi'(\gamma(T-t))
v_\nu}\\[0.3 cm]
\displaystyle{\qquad\qquad\qquad\qquad\qquad\qquad\quad+ \sum_{jk}
\partial_{x_j}([\varphi_\nu,\; a_{jk}]\partial_{x_k} v) \|_{L^2}^2\;
dt}\\[0.3 cm]
\displaystyle{\qquad\geq {1\over K}\int_0^{T\over 2}
\sum_\nu\|\partial_tv_\nu+\sum_{jk}
\partial_{x_j}(a_{jk}\partial_{x_k} v_\nu) + \Phi'(\gamma(T-t))
v_\nu\|_{L^2}^2\;dt}\\[0.5 cm]
\displaystyle{\qquad\qquad\qquad\qquad-{1\over K}\int_0^{T\over
2}\sum_\nu\|\sum_{jk}
\partial_{x_j}([\varphi_\nu,\; a_{jk}]\partial_{x_k} v_\nu) \|_{L^2}^2\; dt}\\
\end{array}
\label{11}
\end{equation} where $[\varphi_\nu,\; a_{jk}] w
=\varphi_\nu(D)(a_{jk}w)-a_{jk}\varphi_\nu(D)w$.
\vskip 0.2 cm
\noindent {\it b) the approximation procedure}
\vskip 0.2 cm
\noindent We start to estimate
$$
\int_0^{T\over 2}\sum_\nu\|\partial_tv_\nu+\sum_{jk}
\partial_{x_j}(a_{jk}\partial_{x_k} v_\nu) + \Phi'(\gamma(T-t))
v_\nu\|_{L^2}^2\;dt.
$$ We obtain
\begin{equation}
\begin{array}{l}
\displaystyle{\int_0^{T\over 2}\sum_\nu\|\partial_tv_\nu+\sum_{jk}
\partial_{x_j}(a_{jk}\partial_{x_k} v_\nu) + \Phi'(\gamma(T-t))
v_\nu\|_{L^2}^2\;dt}\\[0.3 cm]
\displaystyle{\quad=\int_0^{T\over
2}\sum_\nu(\|\partial_tv_\nu\|_{L^2}^2+\|\sum_{jk}
\partial_{x_j}(a_{jk}\partial_{x_k} v_\nu) + \Phi'(\gamma(T-t))
v_\nu\|_{L^2}^2}\\[0.4 cm]
\displaystyle{\qquad+\gamma\Phi''(\gamma(T-t))\|v_\nu\|_{L^2}^2+2\,{\rm
Re}\,\langle
\partial_tv_\nu,\; \sum_{jk} \partial_{x_j}(a_{jk}\partial_{x_k}
v_\nu)\rangle_{L^2})\;dt.}\\
\end{array}
\label{12}
\end{equation} We remark that if $a_{jk}$ would be Lipschitz-continuous the
last term in (\ref{12})
would be easily computed by integration by parts. On the contrary here we
approximate it using a
technique similar to the one of \cite{CDGS} (see also \cite{CL} and
\cite{CDS}). Let $\rho\in
C^\infty_0({\mathbb R})$ with $\,{\rm supp}\, \rho\subseteq[-1/2,\; 1/2]$,
$\int_{\mathbb R} \rho(s)\;ds=1$ and
$\rho(s)\geq 0$ for all $s\in {\mathbb R}$; we set
$$ a_{jk,\,\varepsilon}(t,x)=\int_{\mathbb R} a_{jk}(s,x){1\over
\varepsilon}\rho({t-s\over s})\; ds
$$ for $\varepsilon\in\;]0, 1/2]$. We obtain that there exist $C$, $\tilde
C>0$ such that
\begin{equation} |a_{jk,\,\varepsilon}(t,x)-a_{jk}(t,x)|\leq C\mu(\varepsilon)
\label{13}
\end{equation} and
\begin{equation} |\partial_t a_{jk,\,\varepsilon}(t,x)|\leq \tilde
C\,{\mu(\varepsilon)\over
\varepsilon}
\label{14}
\end{equation}
\vskip0.2 cm
\noindent for all $j,k=1\dots,n$ and for all $(t,x)\in [0,T]\times{\mathbb
R}^n_x$ . We have
$$
\begin{array}{l}
\displaystyle{\int_0^{T\over 2}2\,{\rm Re}\,\langle
\partial_tv_\nu,\; \sum_{jk} \partial_{x_j}(a_{jk}\partial_{x_k}
v_\nu)\rangle_{L^2}\;dt}\\[0.3 cm]
\displaystyle{\qquad=-2\,{\rm Re}\,\int_0^{T\over 2}\sum_{jk}
\langle\partial_{x_j}\partial_t v_\nu,
\; a_{jk}\partial_{x_k} v_\nu\rangle_{L^2}\;dt}\\[0.3 cm]
\displaystyle{\qquad=-2\,{\rm Re}\,\int_0^{T\over 2}\sum_{jk}
\langle\partial_{x_j}\partial_t v_\nu,
\; (a_{jk}-a_{jk,\,\varepsilon})\partial_{x_k}
v_\nu\rangle_{L^2}\;dt}\\[0.3 cm]
\displaystyle{\qquad\qquad\qquad\qquad -2\,{\rm Re}\,\int_0^{T\over 2}\sum_{jk}
\langle\partial_{x_j}\partial_t v_\nu, \; a_{jk,\,\varepsilon}\partial_{x_k}
v_\nu\rangle_{L^2}\;dt.}\\
\end{array}
$$ We remark that $\|\partial_{x_j}v_\nu\|_{L^2}\leq
2^{\nu+1}\|v_\nu\|_{L^2}$ and
$\|\partial_{x_j}\partial_tv_\nu\|_{L^2}\leq
2^{\nu+1}\|\partial_tv_\nu\|_{L^2}$ for all $\nu\in
{\mathbb N}$ so that from (\ref{13}) we get
$$
\begin{array}{l}
\displaystyle{|2\,{\rm Re}\,\int_0^{T\over 2}\sum_{jk}
\langle\partial_{x_j}\partial_t v_\nu, \;
(a_{jk}-a_{jk,\,\varepsilon})\partial_{x_k} v_\nu\rangle_{L^2}\;dt|}\\[0.3 cm]
\displaystyle{\qquad\leq C \mu(\varepsilon)\int_0^{T\over 2}\sum_{jk}
\|\partial_{x_j}\partial_t v_\nu\|_{L^2}\;\|\partial_{x_k}v_\nu\|_{L^2}\;
dt}\\[0.3 cm]
\displaystyle{\qquad\qquad\leq {n^2C\over N}\int_0^{T\over
2}\|\partial_tv_\nu\|_{L^2}^2\;dt  +n^2CN\,
2^{4(\nu+1)}\mu(\varepsilon)\int_0^{T\over 2}\|v_\nu\|_{L^2}^2\;dt} \\
\end{array}
$$ for all $N>0$, and similarly from (\ref{14}) we deduce
$$
\begin{array}{l}
\displaystyle{|2\,{\rm Re}\,\int_0^{T\over 2}
\sum_{jk}\langle\partial_{x_j}\partial_t v_\nu, \;
a_{jk,\,\varepsilon}\partial_{x_k} v_\nu\rangle_{L^2}\;dt|=|\int_0^{T\over
2}\sum_{jk}\langle\partial_{x_j}v_\nu, \;
\partial_t a_{jk,\,\varepsilon}\partial_{x_k} v_\nu\rangle_{L^2}\;dt|}
\\[0.5 cm]
\displaystyle{\qquad\qquad
\leq n \tilde C\, {\mu(\varepsilon)\over \varepsilon}\int_0^{T\over 2}
\|\nabla v_\nu\|_{L^2}^2\; dt\leq n^2 \tilde C\, 2^{2(\nu+1)}\,
{\mu(\varepsilon)\over
\varepsilon}\int_0^{T\over 2}
\|v_\nu\|_{L^2}^2\; dt.}\\
\end{array}
$$ Let $N=n^2C$. We deduce that, for all $\nu\in {\mathbb N}$,
\begin{equation}
\begin{array}{l}
\displaystyle{\int_0^{T\over 2}\|\partial_tv_\nu+\sum_{jk}
\partial_{x_j}(a_{jk}\partial_{x_k} v_\nu) + \Phi'(\gamma(T-t))
v_\nu\|^2_{L^2}\;dt}\\[0.3 cm]
\displaystyle{\qquad\geq\int_0^{T\over 2}(\|\sum_{jk}
\partial_{x_j}(a_{jk}\partial_{x_k} v_\nu) + \Phi'(\gamma(T-t)) v_\nu\|^2_{L^2}
+\gamma\Phi''(\gamma(T-t))\|v_\nu\|^2_{L^2}}\\[0.3 cm]
\displaystyle{\qquad\qquad\qquad\qquad-(n^4C^2\,
2^{4(\nu+1)}\,\mu(\varepsilon)+n^2 \tilde C\,
2^{2(\nu+1)}\, {\mu(\varepsilon)\over \varepsilon})
\|v_\nu\|^2_{L^2})\;dt.}\\
\end{array}
\label{15}
\end{equation} Let $\nu=0$. From (\ref{6bis}) we can choose $\gamma_0>0$
such that
$\Phi''(\gamma(T-t))\geq 1$ for all $\gamma>\gamma_0$ and for all
$t\in [0,\, T/2]$. Taking now $\varepsilon=1/2$ we obtain from (\ref{15}) that
$$
\begin{array}{l}
\displaystyle{\int_0^{T\over 2}\|\partial_tv_0+\sum_{jk}
\partial_{x_j}(a_{jk}\partial_{x_k} v_0) + \Phi'(\gamma(T-t))
v_0\|^2_{L^2}\;dt}\\[0.3 cm]
\displaystyle{\qquad\qquad\qquad\qquad\geq\int_0^{T\over
2}(\gamma-8n^2\mu({1\over 2})(2n^2C^2+\tilde
C))\|v_0\|^2_{L^2}\;dt}\\
\end{array}
$$ for all $\gamma>\gamma_0$. Possibly choosing a larger $\gamma_0$ we
have, again for all
$\gamma>\gamma_0$,
\begin{equation}\begin{array}{l}\displaystyle{
\int_0^{T\over 2}\|\partial_tv_0+\sum_{jk}
\partial_{x_j}(a_{jk}\partial_{x_k} v_0) + \Phi'(\gamma(T-t))
v_0\|^2_{L^2}\;dt}\\[0.3
cm]\displaystyle{\qquad\qquad\qquad\qquad\qquad\qquad\qquad\qquad\qquad
\geq {\gamma\over 2} \int_0^{T\over 2}\|v_0\|^2_{L^2}\; dt.}\\\end{array}
\label {16}
\end{equation} Let now $\nu\geq 1$. We recall that in this case
$\|\nabla v_\nu\|\geq
2^{\nu-1}\|v_\nu\|$. We take $\varepsilon=2^{-2\nu}$.  We obtain from
(\ref{15}) that
\begin{equation}
\begin{array}{l}
\displaystyle{\int_0^{T\over 2}\|\partial_tv_\nu+\sum_{jk}
\partial_{x_j}(a_{jk}\partial_{x_k} v_\nu) + \Phi'(\gamma(T-t))
v_\nu\|^2_{L^2}\;dt}\\[0.3 cm]
\displaystyle{\qquad\geq\int_0^{T\over 2}(\|\sum_{jk}
\partial_{x_j}(a_{jk}\partial_{x_k} v_\nu) + \Phi'(\gamma(T-t)) v_\nu\|^2_{L^2}
}\\[0.3 cm]
\displaystyle{\qquad\qquad\qquad+\gamma\Phi''(\gamma(T-t))\|v_\nu\|^2_{L^2}-K\,
2^{4\nu}\,\mu(2^{-2\nu})\|v_\nu\|^2_{L^2})\; dt}\\[0.3 cm]
\displaystyle{\qquad\geq\int_0^{T\over 2}((\|\sum_{jk}
\partial_{x_j}(a_{jk}\partial_{x_k} v_\nu)\|_{L^2} - \Phi'(\gamma(T-t))\|
v_\nu\|_{L^2})^2 }\\[0.3 cm]
\displaystyle{\qquad\qquad\qquad+\gamma\Phi''(\gamma(T-t))\|v_\nu\|^2_{L^2}-K\,
2^{4\nu}\,\mu(2^{-2\nu})\|v_\nu\|^2_{L^2})\; dt}
\end{array}
\label{17}
\end{equation} where $K=16n^4C^2+4n^2\tilde C$.  On the other hand we have
$$
\begin{array}{l}
\displaystyle{ \|\sum_{jk} \partial_{x_j}(a_{jk}\partial_{x_k} v_\nu)\|_{L^2}
\;\|v_\nu\|_{L^2}\geq
|\langle\sum_{jk} \partial_{x_j}(a_{jk}\partial_{x_k} v_\nu),\;
v_\nu\rangle_{L^2}|}\\[0.3 cm]
\displaystyle{\qquad\qquad\geq|\sum_{jk}\langle a_{jk}\partial_{x_k} v_\nu,\;
\partial_{x_j}v_\nu\rangle_{L^2}|
\geq \lambda_0\|\nabla v_\nu\|^2_{L^2}\geq {\lambda_0\over 4}\,
2^{2\nu}\,\| v_\nu\|^2_{L^2}}\\
\end{array}
$$ and consequently
\begin{equation}
\|\sum_{jk} \partial_{x_j}(a_{jk}\partial_{x_k} v_\nu)\|_{L^2} \geq
{\lambda_0\over 4}\, 2^{2\nu}\,\|_{L^2}
v_\nu\|.
\label{18}
\end{equation} Suppose first that $\Phi'(\gamma(T-t))\leq{\lambda\over 8}\,
2^{2\nu}$. Then from
(\ref{18}) we deduce that
$$
 \|\sum_{jk} \partial_{x_j}(a_{jk}\partial_{x_k} v_\nu)\|_{L^2} -
\Phi'(\gamma(T-t))\| v_\nu\|_{L^2}\geq
{\lambda\over 8}\, 2^{2\nu}\|v_\nu\|_{L^2}
$$ and then, using also the fact that $\Phi''(\gamma(T-t))\geq 1$, we obtain
that there exist
$\gamma_0$ and  $c>0$ such that\begin{equation}
\begin{array}{l}
\displaystyle{\int_0^{T\over 2}((\|\sum_{jk}
\partial_{x_j}(a_{jk}\partial_{x_k} v_\nu)\|_{L^2} - \Phi'(\gamma(T-t))\|
v_\nu\|_{L^2})^2 }\\[0.3 cm]
\displaystyle{\qquad\qquad\qquad+\gamma\Phi''(\gamma(T-t))\|v_\nu\|_{L^2}^2-K\,
2^{4\nu}\,\mu(2^{-2\nu})\|v_\nu\|_{L^2}^2\; )dt}\\[0.3 cm]
\displaystyle{\quad \geq \int_0^{T\over 2}(({\lambda\over 8}\,
2^{2\nu})^2+\gamma -K\,
2^{4\nu}\,\mu(2^{-2\nu})\|v_\nu\|_{L^2}^2)\; dt}\\[0.5 cm]
\displaystyle{\qquad\geq\int_0^{T\over 2}(({\lambda\over 16})^2\,
2^{4\nu}+{2\over 3}\gamma)\|v_\nu\|_{L^2}^2\; dt \geq
\int_0^{T\over 2}({\gamma\over 2}+c\gamma^{1\over 2}\,
2^{2\nu})\|v_\nu\|_{L^2}^2\; dt}\\
\end{array}
\label{19}
\end{equation} for all $\gamma\geq \gamma_0$. If on the contrary
$\Phi'(\gamma(T-t))\geq{\lambda\over 8}\, 2^{2\nu}$ then, using (\ref{6}),
the fact that
$\lambda_0\leq 1$ and the properties of $\mu$ (see Remark \ref{rem1}),
$$
\begin{array}{l}
\displaystyle{\Phi''(\gamma(T-t))= (\Phi'(\gamma(T-t))^2\mu({1\over
\Phi'(\gamma(T-t))})}\\[0.3 cm]
\displaystyle{\qquad\qquad\qquad\geq ({\lambda_0\over 8})^2\, 2^{4\nu}\,
\mu({8\over
\lambda_0}\,2^{-2\nu})
\geq ({\lambda_0\over 8})^2\, 2^{4\nu}\, \mu(2^{-2\nu}).}\\
\end{array}
$$
Hence also in this case there exist $\gamma_0$ and  $c>0$ such that
\begin{equation}
\begin{array}{l}
\displaystyle{\int_0^{T\over 2}((\|\sum_{jk}
\partial_{x_j}(a_{jk}\partial_{x_k} v_\nu)\|_{L^2} - \Phi'(\gamma(T-t))\|
v_\nu\|_{L^2})^2 }\\[0.3 cm]
\displaystyle{\qquad\qquad\qquad+\gamma\Phi''(\gamma(T-t))\|v_\nu\|_{L^2}^2-K\,
2^{4\nu}\,\mu(2^{-2\nu})\|v_\nu\|_{L^2}^2)\; dt}\\[0.3 cm]
\displaystyle{\quad \geq \int_0^{T\over 2}({\gamma\over 2}+({\gamma\over
2}({\lambda\over
8})^2-K)\, 2^{4\nu}\mu(2^{-2\nu}))\|v_\nu\|_{L^2}^2\; dt}\\[0.3cm]
\displaystyle{\qquad\geq \int_0^{T\over 2}({\gamma\over 2}+c\gamma\,
2^{2\nu})\|v_\nu\|_{L^2}^2\; dt}\\
\end{array}
\label{20}
\end{equation} for all $\gamma\geq \gamma_0$. Putting together (\ref{19})
and (\ref{20}) we have
that there exist $\gamma_0$ and $c>0$ such that
\begin{equation}
\begin{array}{l}
\displaystyle{\int_0^{T\over 2}\|\partial_t v_\nu+ \sum_{jk}
\partial_{x_j}(a_{jk}\partial_{x_k} v_\nu) + \Phi'(\gamma(T-t)) v_\nu\|_{L^2}^2
\, dt}\\[0.3 cm]
 \displaystyle{\qquad \geq\int_0^{T\over 2}({\gamma\over 2}+c\gamma^{1\over
2}\,
2^{2\nu})\|v_\nu\|_{L^2}^2\; dt}\\
\end{array}
\label{21}
\end{equation} for all $\nu\geq 1$ and for all $\gamma\geq \gamma_0$.

Form (\ref{16}) and (\ref{21}) we get that there exist $\gamma_0$ and
$\tilde c>0$ such that
\begin{equation}
\begin{array}{l}
\displaystyle{\int_0^{T\over 2}\sum_\nu\|\partial_t v_\nu+ \sum_{jk}
\partial_{x_j}(a_{jk}\partial_{x_k} v_\nu) + \Phi'(\gamma(T-t)) v_\nu\|_{L^2}^2
\, dt}\\[0.3 cm]
 \displaystyle{\qquad \geq\tilde c \gamma^{1\over 2}\int_0^{T\over
2}\sum_\nu(\gamma^{1\over
2}\|v_\nu\|_{L^2}^2 +
\|\nabla v_\nu\|_{L^2}^2)\; dt}\\
\end{array}
\label{22}
\end{equation} for all $\gamma\geq \gamma_0$.

\vskip 0.2 cm

\noindent {\it c) the estimate for the commutator}
\vskip 0.2 cm
\noindent
\noindent For $\psi\in C^\infty_0({\mathbb R}^n_\xi)$, we define $\check
\psi (x) =
{1\over
(2\pi)^n}\int_{{\mathbb R}^n_\xi} e^{{\rm i}x\xi}\psi(\xi)\; d\xi$. Notice that
$(\nabla\psi)\! \check{\phantom
I}\!(x)=i\check \psi(x) x$ and
${(\psi_1\psi_2)\! \check{\phantom I}\!} =\check \psi_1\ast\check \psi_2$.
For $w\in L^2({\mathbb R}^n_x,{\mathbb C})$ we have
$$ w_\nu(x) = \int_{{\mathbb R}^n_y} \check \varphi_\nu (x-y)w(y)\; dy.
$$  Moreover
$$ [\varphi_\nu, \; a_{jk}] w(x)=\int_{{\mathbb R}^n_y} h_{jk}^\nu(x,y)
w(y)\; dy,
$$ where
$$ h_{jk}^\nu(x,y)= \check \varphi_\nu (x-y)(a_{jk}(y)-a_{jk}(x))
$$ (to avoid cumbersome  notations here and throughout this point we drop
$t$ in writing the
variables of the coefficients $a_{jk}$). One can rewrite $h_{jk}^\nu$ as
$h_{jk}^\nu(x,y)=h_{jk}^{\nu, 1}(x,y)+h_{jk}^{\nu, 2}(x,y)$, where
$$
\begin{array}{l}
\displaystyle{h_{jk}^{\nu, 1}(x,y)=\check \varphi_\nu (x-y)\int_0^1(\nabla
a_{jk}(x+\theta(y-x))-\nabla a_{jk}(x))\cdot(y-x)\; d\theta}\\[0.3 cm]
\displaystyle{h_{jk}^{\nu, 2}(x,y)=\check \varphi_\nu (x-y)\nabla
a_{jk}(x)\cdot(y-x).}\\
\end{array}
$$ We remark that
$$
\int_{{\mathbb R}^n_y}h_{jk}^{\nu, 2}(x,y)w(y)\;
dy=\sum_{\mu=0}^{+\infty}\int_{{\mathbb R}^n_y} h_{jk}^{\nu,
2}(x,y)w_\mu(y)\; dy,
$$ where $w_\mu (x)= \varphi_\mu(D)w(x)$. We have then
$$
\begin{array}{l}
\displaystyle{\int_{{\mathbb R}^n_y} h_{jk}^{\nu, 2}(x,y)w_\mu(y)\;
dy}\\[0.3 cm]
\qquad\displaystyle{=\int_{{\mathbb R}^n_y} \check \varphi_\nu (x-y)\nabla
a_{jk}(x)\cdot(y-x)(\int_{{\mathbb R}^n_z}
\check \varphi_\mu (y-z)w(z)\; dz)\; dy}\\[0.3 cm]
\qquad\displaystyle{=\int_{{\mathbb R}^n_z}\nabla a_{jk}(x)\cdot(\int
_{{\mathbb R}^n_y} \check \varphi_\mu
(y-z)
\check \varphi_\nu (x-y)\; dy)w(z)\; dz}\\[0.3 cm]
\qquad\displaystyle{=\int_{{\mathbb R}^n_z}\nabla
a_{jk}(x)\cdot(\int_{{\mathbb R}^n_y}
{\rm i}  \check \varphi_\mu
(y-z) (\nabla \varphi_\nu)\! \check{\phantom I}\!(x-y)\; dy) w(z) \;
dz}\\[0.3 cm]
\qquad\displaystyle{=\int_{{\mathbb R}^n_z}\nabla a_{jk}(x)\cdot
(\int_{{\mathbb R}^n_y} {\rm i}  \check \varphi_\mu
(y)  (\nabla \varphi_\nu)\! \check{\phantom I}\!((x-z)-y)\; dy) w(z) \; dz.}\\
\end{array}
$$ Recalling that if $\mu<\nu-1$ or $\mu>\nu+1$ then
$$
\int_{{\mathbb R}^n_y}   \check \varphi_\mu (y)  (\nabla \varphi_\nu)\!
\check{\phantom I}\!((x-z)-y)\; dy=
(\varphi_\mu \nabla\varphi_\nu)\! \check{\phantom I}\!(x-z)=0,
$$ we finally obtain
$$
\int_{{\mathbb R}^n_y} h_{jk}^{\nu, 2}(x,y)w(y)\; dy= \int_{{\mathbb R}^n_y}
h_{jk}^{\nu, 2}(x,y)(w_{\nu-1}(y)+
w_\nu(y)+w_{\nu+1}(y))\; dy,
$$ where we have set $w_{-1}= 0$ identically. We deduce
\begin{equation}
\begin{array}{l}
\displaystyle{\partial_{x_l}\, [\varphi_\nu\; a_{jk}] w(x)}\\[0.3 cm]
\qquad\qquad\displaystyle{=\int_{{\mathbb R}^n_y} \partial_{x_l} h_{jk}^{\nu,
1}(x,y)w(y)\; dy }\\[0.3 cm]
\qquad\qquad\qquad\displaystyle{+
\int_{{\mathbb R}^n_y} \partial_{x_l} h_{jk}^{\nu, 2}(x,y)(w_{\nu-1}(y)+
w_\nu(y)+w_{\nu+1}(y))\; dy.}\\
\end{array}
\label{2.1}
\end{equation} Using the explicit expression of $h_{jk}^{\nu, 1}$ we get
$$
\begin{array}{l}
\displaystyle{\partial_{x_l}h_{jk}^{\nu, 1}(x,y)}\\[0.3 cm]
\displaystyle{=\partial_{x_l}\check \varphi_\nu(x-y)\int_0^1(\nabla
a_{jk}(x+\theta(y-x))-\nabla
a_{jk}(x))\cdot (y-x)\; d\theta}\\[0.3 cm]
\quad\displaystyle{+\check
\varphi_\nu(x-y)\int_0^1((1-\theta)\nabla(\partial_{x_l}
a_{jk})(x+\theta(y-x))-\nabla(\partial_{x_l} a_{jk})(x))\cdot (y-x)\,
d\theta}\\[0.3 cm]
\qquad\displaystyle{-\check \varphi_\nu(x-y)\int_0^1 (\partial_{x_l}
a_{jk}(x+\theta(y-x))-\partial_{x_l} a_{jk}(x))\; d\theta.}\\
\end{array}
$$ Using the mean value theorem we deduce that
$$ |\partial_{x_l}h_{jk}^{\nu, 1}(x,y)|\leq (|\partial_{x_l}\check
\varphi_\nu(x-y)||x-y|^2+3
|\check
\varphi_\nu(x-y)||x-y|) \|D^2 a_{jk}\|_{L^\infty} .
$$ Hence both $\int_{{\mathbb R}^n_x}|\partial_{x_l}h_{jk}^{\nu, 1}(x,y)|\,
dx$ and
$\int_{{\mathbb R}^n_y}|\partial_{x_l}h_{jk}^{\nu,
1}(x,y)|\, dy$ are dominated by the quantity
\begin{equation}
\|D^2 a_{jk}\|_{L^\infty}  \int_{{\mathbb R}^n}(|\partial_{x_l}\check
\varphi_\nu(z)||z|^2+3 |\check
\varphi_\nu(z)||z|)\, dz.
\label {2.2}
\end{equation} Now we observe that, for $\nu\geq 1$,
\begin{equation}
\check \varphi_\nu(z)=\check \varphi( 2^\nu z)\, 2^{n\nu}\quad {\rm
and}\quad \partial_{x_l}\check
\varphi_\nu(z)=
 \partial_{x_l}\check \varphi(2^\nu z)\, 2^{(n+1)\nu},
\label{2.3}
\end{equation} where $\varphi(\xi)=\varphi_0(\xi)-\varphi_0(2\xi)$. Setting
$2^\nu z=\zeta$, the
quantity in (\ref{2.2}) becomes
$$ 2^{-\nu}\|D^2 a_{jk}\|_{L^\infty}  \int_{{\mathbb
R}^n_\zeta}(|\partial_{x_l}\check
\varphi(\zeta)||\zeta|^2+ 3
|\check
\varphi(\zeta)||\zeta|)\, d\zeta.
$$ Consequently there exists $K>0$ such that, for all $\nu\geq 0$,
\begin{equation}
\|\int_{{\mathbb R}^n_y}\partial_{x_l}h_{jk}^{\nu, 1}(\cdot,y)w(y)\;
dy\|_{L^2}\leq K2^{-\nu} \|w\|_{L^2}.
\label{2.4}
\end{equation} Next we consider
$$
\begin{array}{l}
\displaystyle{\partial_{x_l}h_{jk}^{\nu,2}(x,y)}= \partial_{x_l}\check
\varphi_\nu(x-y)
\nabla a_{jk}(x)\cdot (y-x)\\[0.3 cm]
\qquad\qquad\qquad\quad\displaystyle{+\check
\varphi_\nu(x-y)\nabla(\partial_{x_l}a_{jk})(x)\cdot
(y-x)}\\[0.3 cm]
\qquad\qquad\qquad\qquad\displaystyle{+\check
\varphi_\nu(x-y)\partial_{x_l}a_{jk}(x).}\\
\end{array}
$$ Again both $\int_{{\mathbb R}^n_x}|\partial_{x_l}h_{jk}^{\nu, 2}(x,y)|\,
dx$ and
$\int_{{\mathbb R}^n_y}|\partial_{x_l}h_{jk}^{\nu,
2}(x,y)|\, dy$ are dominated by
\begin{equation}
\begin{array}{l}
\displaystyle{\|\nabla a_{jk}\|_{L^\infty}\int_{{\mathbb R}^n_z}|
\partial_{x_l}\check
\varphi_\nu(z)||z|\, dz+
\| D^2 a_{jk}\|_{L^\infty}\int_{{\mathbb R}^n_z}| \check \varphi_\nu(z)||z|\,
dz}\\[0.3 cm]
\qquad\displaystyle{+\|\nabla a_{jk}\|_{L^\infty}\int_{{\mathbb R}^n_z}| \check
\varphi_\nu(z)||z|\, dz.}\\
\end{array}
\label {2.5}
\end{equation} As before setting $2^\nu z=\zeta$ and recalling (\ref{2.3})
we have that there exists
$K>0$ such that, for all $\nu\geq 0$,
\begin{equation}
\|\int_{{\mathbb R}^n_y}\partial_{x_l}h_{jk}^{\nu, 2}(\cdot,y)w(y)\;
dy\|_{L^2}\leq K \|w\|_{L^2}.
\label{2.6}
\end{equation} It follows from (\ref{2.1}), (\ref{2.4}) and (\ref{2.6}) that
$$
\|\partial_{x_l} [\varphi_\nu,\, a_{jk}] w\|_{L^2}\leq
K(2^{-\nu}\,\|w\|_{L^2}+\|w_{\nu-1}\|_{L^2}+\|w_\nu\|_{L^2}+\|w_{\nu+1}\|_{L^2})
$$ for all $\nu\geq 0$.  Hence, possibly choosing a larger $K>0$,
$$
\begin{array}{l}
\displaystyle{\|\partial_{x_j} [\varphi_\nu,\,
a_{jk}]\partial_{x_k}v\|_{L^2}\leq
K(2^{-\nu}\|\partial_{x_k}v\|_{L^2} +\|(\partial_{x_k}v)_{\nu-1}\|_{L^2}}
\\[0.3 cm]
\qquad\qquad\qquad\qquad\qquad\qquad\qquad\displaystyle{
+\|(\partial_{x_k}v)_{\nu}\|_{L^2}+\|(\partial_{x_k}v)_{\nu+1}\|_{L^2}) }\\
\end{array}
$$ for all $j$, $k=1,\dots,n$,  $\nu\geq 0$ and $v\in C^\infty_0({\mathbb
R}^n,{\mathbb C})$. Finally from
(\ref{10}) we obtain that there exists a $\tilde K$ such that
\begin{equation}
\sum_\nu\|\sum_{jk}\partial_{x_j} [\varphi_\nu,\,
a_{jk}]\partial_{x_k}v\|_{L^2}^2\leq \tilde
K\|\nabla v\|^2_{L^2}.
\label{2.7}
\end{equation}
\vskip 0.2 cm \noindent {\it d) end of the proof of Proposition 1}
\vskip 0.2 cm \noindent From (\ref{11}), (\ref{22}) and
(\ref{2.7}) we obtain that there exist $\gamma_0$, $\tilde c$, $K$
and $\tilde K$ positive constants such that $$
\begin{array}{ll}
\displaystyle{\int_0^{T\over 2} \|\partial_tv+\sum_{jk}
\partial_{x_j}(a_{jk}\partial_{x_k} v) + \Phi'(\gamma(T-t))
v\|^2_{L^2}\;
dt}\\[0.3 cm]
\displaystyle{\qquad\geq {\tilde c\over K}\gamma^{1\over 2}\int_0^{T\over 2}
\sum _\nu (\|\nabla v_\nu\|^2_{L^2}+\gamma^{1\over
2}\|v_\nu\|^2_{L^2})\; dt- {\tilde K\over K}\int_0^{T\over 2}\|\nabla
v\|^2_{L^2}\; dt}
\end{array}$$for all $v\in C^\infty_0({\mathbb R}_t \times {\mathbb R}^n_x,
{\mathbb C})$ with support
in $[0,{T/ 2}]\times {\mathbb R}^n_x$ and for all $\gamma\geq \gamma_0$.
Using (\ref{10}) we immediately obtain (\ref{8}) and
the proof of the Proposition \ref{p1} is complete.

Let us come finally to the proof of  Theorem \ref{t1}. First of all we remark
that a density argument ensures that the inequality (\ref{7}) holds for all
$\gamma\geq \gamma_0$ and for all $u\in {\mathcal
H}_1$ such that ${\rm supp}\, u\subseteq [0,{T/ 2}]\times {\mathbb R}^n_x$.
Suppose now that $u\in {\mathcal H}_1$,
$u(0,x)=0$ in ${\mathbb R}^n_x$ and \begin{equation}\|\partial_t u +
\sum_{jk}\partial_{x_j}(a_{jk}\partial_{x_k}
u)\|^2_{L^2}\leq\tilde C(  \|\nabla u\|^2_{L^2}+\|u\|^2_{L^2})\label
{2.8}\end{equation}for a.e. $t\in [0,T]$. We consider
$\omega\in C^\infty ({\mathbb R}_t)$ such that $\omega(t)=0$ for all $t\geq
T/2$ and $\omega(t)=1$ for all $t\leq T/3$. We
apply (\ref{7}) to the function $\omega(t)u(t,x)$ and we obtain
$$
\begin{array}{l}
\displaystyle{\int_0^{T\over 2} e^{{2\over \gamma}\Phi(\gamma(T-t))}
\| \partial_t(\omega u)+\sum_{jk}\partial_{x_j}(a_{jk}\partial_{x_k}
(\omega u))\|^2_{L^2}\; dt}\\[0.3
cm]\displaystyle{\qquad\qquad\qquad \geq C\gamma^{1\over 2}\int_0^{T\over
2} e^{{2\over \gamma}\Phi(\gamma(T-t))}(\|\nabla
(\omega u)\|^2_{L^2}+\gamma^{1\over 2}\|\omega u\|^2_{L^2})\; dt}\\
\end{array}
$$
and consequently$$\begin{array}{l}
\displaystyle{\int_0^{T\over 3} e^{{2\over \gamma}\Phi(\gamma(T-t))}
\| \partial_t u+\sum_{jk}\partial_{x_j}(a_{jk}\partial_{x_k} u)\|^2_{L^2}\;
dt}\\[0.3 cm]\displaystyle{\quad +\int_{T\over
3}^{T\over 2} e^{{2\over \gamma}\Phi(\gamma(T-t))}\| \partial_t(\omega
u)+\sum_{jk}\partial_{x_j}(a_{jk}\partial_{x_k}
(\omega u))\|^2_{L^2}\; dt}\\[0.3 cm]\displaystyle{\qquad\qquad\qquad\qquad
\geq C\gamma^{1\over 2}\int_0^{T\over 3}
e^{{2\over \gamma}\Phi(\gamma(T-t))}(\|\nabla u\|^2_{L^2}+\gamma^{1\over
2}\| u\|^2_{L^2})\; dt.}\\
\end{array}
$$
By
(\ref{2.8}) we get
$$
\begin{array}{l}
\displaystyle{\int_{T\over 3}^{T\over 2} e^{{2\over \gamma}
\Phi(\gamma(T-t))}\| \partial_t(\omega
u)+\sum_{jk}\partial_{x_j}(a_{jk}\partial_{x_k} (\omega u))\|^2_{L^2}\;
dt}\\[0.3
cm]\displaystyle{\qquad\qquad\qquad \geq \int_0^{T\over 3} e^{{2\over
\gamma}\Phi(\gamma(T-t))}((C\gamma^{1\over 2}-\tilde
C)\|\nabla u\|^2_{L^2}+(C\gamma-\tilde C)\| u\|^2_{L^2})\; dt,}\\
\end{array}
$$ so that, since $\Phi$ is
increasing,
$$
\begin{array}{l}
\displaystyle{e^{{2\over \gamma}\Phi({2\over 3}\gamma T)}\int_{T\over
3}^{T\over 2}
\| \partial_t(\omega u)+\sum_{jk}\partial_{x_j}(a_{jk}\partial_{x_k}
(\omega u))\|^2_{L^2}\; dt}\\[0.3
cm]\displaystyle{\qquad\qquad\qquad \geq e^{{2\over \gamma}\Phi({3\over
4}\gamma T)}\int_0^{T\over 4} ((C\gamma^{1\over
2}-\tilde C)\|\nabla u\|^2_{L^2}+(C\gamma-\tilde C)\| u\|^2_{L^2})\; dt.}\\
\end{array}
$$ Choosing $\gamma_0$ sufficiently
large we deduce that for all $\gamma\geq \gamma_0$,
$$
\begin{array}{l}
\displaystyle{\int_{T\over 3}^{T\over 2} \| \partial_t(\omega u)+
\sum_{jk}\partial_{x_j}(a_{jk}\partial_{x_k} (\omega u))\|^2_{L^2}\; dt}\\[0.3
cm]\displaystyle{\qquad\qquad\qquad\qquad\qquad \geq {C\over 2}\gamma
e^{{2\over \gamma}(\Phi({3\over 4}\gamma
T)-\Phi({2\over 3}\gamma T))}\int_0^{T\over 4}  \| u\|^2_{L^2}\;
dt.}\\\end{array}
$$
Remarking now
that
$$\lim_{\gamma\to +\infty} {2\over \gamma}(\Phi({3\over 4}\gamma
T)-\Phi({2\over 3}\gamma T))=\lim_{\gamma\to
+\infty}{2\over \gamma}\int_{{2\over 3}\gamma T}^{{3\over 4}\gamma
T}\phi^{-1}(\tau)\; d\tau=+\infty,$$we let $\gamma$ go to
$+\infty$ and we deduce that $u(t,x)=0$ in $[0,T/4]\times {\mathbb R}^n_x$.
The conclusion of the proof of the Theorem
\ref{t1} easily follows.

To prove Theorem \ref{t2} it will be sufficient to multiply $u$ by a
function $\theta\in
C^\infty({\mathbb R}^n_x)$ such that $\theta>0$ and $\theta(x)=e^{-2C|x|}$
for all $x\in {\mathbb R}^n_x$  with $|x|\geq 1$.
A direct computation shows that $\theta u\in {\mathcal H}_1$ and satisfies
(\ref{2.8}). Consequently $\theta u=0$ in
$[0,T]\times {\mathbb R}^n_x$ and the same will be for $u$.

\section {Sketch of the proof of Theorem 3}

In the proof of Theorem \ref{t3} we will follow closely the
construction of the example in \cite {Pl}. Let $A$, $B$, $C$, $J$
be four $C^\infty$ functions defined in $\mathbb R$ with $0\leq
A(s),\ B(s),\ C(s)\leq 1$, $-2\leq J(s)\leq 2$ for all $s\in
{\mathbb R}$ and $$
\begin{array}{ll}
\displaystyle{A(s)=1\quad{\rm for}\ s\leq {1\over 5}, }&
\quad\displaystyle{A(s)=0\quad{\rm for}\
s\geq {1\over 4},}\\[0.3 cm]
\displaystyle{B(s)=0\quad{\rm for}\ s\leq 0\ {\rm or}\ s\geq 1, }&\quad
\displaystyle{B(s)=1\quad{\rm for}\ {1\over 6}\leq s\leq {1\over 2},}\\[0.3 cm]
\displaystyle{C(s)=0\quad{\rm for}\ s\leq {1\over 4}, }&
\quad\displaystyle{C(s)=1\quad{\rm for}\
s\geq {1\over 3},}\\[0.3 cm]
\displaystyle{J(s)=-2\quad{\rm for}\ s\leq {1\over 6}\ {\rm or}\ s\geq
{1\over 2}, }&\quad
\displaystyle{J(s)=2\quad{\rm for}\ {1\over 5}\leq s\leq {1\over 3}.}\\
\end{array}
$$ Let $(a_n)_n$, $(z_n)_n$ be two real sequences such that
\begin{eqnarray}
\displaystyle{-1<a_n<a_{n+1}\quad {\rm for\ all}\ n\geq
1},&&\displaystyle{\lim_n a_n=0,}\label
{4.1}\\[0.3 cm]
\displaystyle{1<z_n<z_{n+1}\quad {\rm for\ all}\ n\geq
1},&&\displaystyle{\lim_n
z_n=+\infty;}\label {4.2}
\end{eqnarray} and let us define $r_n=a_{n+1}-a_n$, $q_1=0$,
$q_n=\sum_{k=2}^nz_k r_{k-1}$ for all
$n\geq 2$, and $p_n=(z_{n+1}-z_n)r_n$. We suppose moreover that
\begin{equation} p_n>1\quad {\rm for \  all}\ n\geq 1.
\label{4.3}
\end{equation} We set $A_n(t)=A({t-a_n\over r_n})$, $B_n(t)=B({t-a_n\over
r_n})$,
$C_n(t)=C({t-a_n\over r_n})$ and $J_n(t)=J({t-a_n\over r_n})$. We define
$$
\begin{array}{l}
\displaystyle{v_n(t,x_1)=\exp (-q_n-z_n(t-a_n))\cos \sqrt {z_n} x_1,}\\[0.3 cm]
\displaystyle{w_n(t,x_2)=\exp (-q_n-z_n(t-a_n)+J_n(t)p_n)\cos \sqrt {z_n}
x_2,}\\
\end{array}
$$ and
$$
\begin{array}{l} u(t,x_1, x_2)\\[0.3 cm] =\left\{
\begin{array}{ll}
\displaystyle{v_1(t, x_1)} &\displaystyle{{\rm for }\ t\leq a_1, }\\[0.3 cm]
\displaystyle{A_n(t)v_n(t,x_1)+B_n(t)w_n(t,x_2)+C_n(t)v_{n+1}(t,x_1)}&\displaystyle{{\rm
for }\
a_n\leq t\leq a_{n+1} ,}\\[0.3 cm]
\displaystyle{0} &\displaystyle{{\rm for }\ t\geq 0. }\\
\end{array}
\right.
\end{array}
$$ If for all $\alpha$, $\beta$ $\gamma>0$
\begin{equation}
\lim_n\exp(-q_n+2p_n)z_{n+1}^\alpha p_n^\beta r_n^{-\gamma}=0
\label{4.4}
\end{equation} then $u$ is a $C^\infty_b({\mathbb R}^3)$ function. We define
$$ l(t)=\left\{
\begin{array}{ll}
\displaystyle{1} &\displaystyle{{\rm for }\ t\leq a_1\ {\rm or }\ t\geq 0,
}\\[0.3 cm]
\displaystyle{1+J'_n(t)p_nz^{-1}_n}&\displaystyle{{\rm for }\ a_n\leq t\leq
a_{n+1} .}\\
\end{array}
\right.
$$ The condition
\begin{equation}
\sup_n\; \{p_nr_n^{-1}z_n^{-1}\}\leq {1\over 2\|J'\|_{L^\infty}}
\label{4.5}
\end{equation} guarantees that the operator $\text{\it
\L}=\partial_t-\partial^2_{x_1}-l(t)\partial^2_{x_2}$ is parabolic.
Moreover $l$ is a $C^\mu$
function under the condition
\begin{equation}
\sup_n\; \{{p_nr_n^{-1}z_n^{-1}\over \mu(r_n)}\}<+\infty.
\label{4.6}
\end{equation} Finally we define
$$
\begin{array}{l}
\displaystyle{b_1=- {\text{\it \L} u\over
u^2+(\partial_{x_1}u)^2+(\partial_{x_2}u)^2}\partial_{x_1}u,}\\[0.3cm]
\displaystyle{b_2=- {\text{\it \L} u\over
u^2+(\partial_{x_1}u)^2+(\partial_{x_2}u)^2}\partial_{x_2}u,}\\[0.3cm]
\displaystyle{c=- {\text{\it \L} u\over
u^2+(\partial_{x_1}u)^2+(\partial_{x_2}u)^2}u}\\
\end{array}
$$ and as in \cite {Pl} the coefficients $b_1$,  $b_2$,  $c$ wil be in
$C^\infty_b$ if for all
$\alpha$, $\beta$, $\gamma>0$
\begin{equation}
\lim_n\exp(-p_n)z_{n+1}^\alpha p_n^\beta r_n^{-\gamma}=0.
\label{4.7}
\end{equation} We choose
\begin{equation} a_n= -\sum _{j=n}^{+\infty}{1\over (j+k_0)^2\mu({1\over
j+k_0})}, \qquad
z_n=(n+k_0)^3
\label{4.8}
\end{equation} with $k_0$ sufficiently large. Changing $t$ in $-t$ the
proof of Theorem \ref{t3}
will be complete as soon as under the choice (\ref{4.8}) the conditions
(\ref{4.1}),...,
(\ref{4.7}) hold. Let's verify this.  Since the function $\sigma \mapsto
1/(\sigma^2\mu
(1/\sigma))$ is decreasing on $[1,+\infty[$ (see Remark \ref{rem1}) we have that
the hypothesis
(\ref{4}) is equivalent to the convergence of the series $\sum_n
((n+k_0)^2\mu({1\over
n+k_0}))^{-1}$ and (\ref{4.1}) follows. Condition (\ref{4.2}) is obvious.
We have, for $n\geq 2$,
$$ q_n=\sum_{j=2}^n(j+k_0)^3{1\over (j+k_0-1)^2\mu({1\over  j+k_0-1})}\geq
\sum_{j=2}^n{j+k_0\over
\mu({1\over  j+k_0-1})}.
$$ Remarking that $\mu({1\over  j+k_0-1})\leq 1$ we obtain that
\begin{equation} q_n\geq {1\over 2} ((n+k_0+1)(n+k_0)-(k_0+3)(k_0+2))
\label{4.9}
\end{equation} for all $n\geq 2$. On the other hand
$$ p_n= ( 3(n+k_0)^2+3(n+k_0)+1){1\over (n+k_0)^2\mu({1\over  n+k_0})};
$$ using also the fact that there exists $c>0$ such that $\mu(s)\geq cs$
for all $s\in [0,1]$ we
deduce that
$$ {3\over \mu({1\over  n+k_0})}\leq p_n\leq {3\over c}(n+k_0+2)
$$ for all $n\geq 1$. Finally remarking that it is not restrictive to
suppose that $\mu(s)\leq
s^{1/2}$ for all $s\in [0,1]$  (if it is not so it is sufficient to replace
$\mu(s)$ with $\min
\;\{\mu(s),\; s^{1\over 2}\}$), we have
\begin{equation} 3(n+k_0)^{1\over 2}\leq p_n\leq {3\over c}(n+k_0+2)
\label {4.10}
\end{equation} and
\begin{equation} (n+k_0)^{-{3\over 2}}\leq r_n\leq {1\over c}(n+k_0)^{-1}
\label {4.11}
\end{equation} for all $n\geq 1$.  Easily the first part of (\ref{4.10})
implies (\ref{4.3}) if
$k_0$ is sufficiently large and  (\ref{4.9}), (\ref{4.10}) and (\ref{4.11})
give (\ref{4.4}) and
(\ref{4.7}). We observe that
\begin{equation}
\begin{array}{l}
\displaystyle{ p_n r_n^{-1} z_n^{-1}=(z_{n+1}-z_n)z_n^{-1}} \\[0.3 cm]
\displaystyle{\qquad=3(n+k_0)^{-1}+3(n+k_0)^{-2} +(n+k_0)^{-3}\leq 7
(n+k_0)^{-1}}\\
\end{array}
\label{4.12}
\end{equation} for all $n\geq 1$ and again taking $k_0$ is sufficiently
large (\ref{4.5}) follows.
To prove (\ref{4.6}) we start remarking that since the function
$s \mapsto {\mu(s )\over s}$ is decreasing on $\;]0,1]$ and $\lim_{s \to
0^+}{\mu(s)\over
s}=+\infty$ (see Remark \ref{rem1}) we have that there exists $k_0$ such that
$$ r_n(n+k_0)={1\over (n+k_0) \mu({1\over n+k_0})}={{1\over (n+k_0)}\over
\mu({1\over n+k_0})}\leq 1
$$
 for all $n\geq 1$, so that $r_n\leq {1\over n+k_0}$ and then
\begin{equation} {\mu(r_n)\over r_n}\geq {\mu({1\over n+k_0})\over {1\over
(n+k_0)}}
\label{4.13}
\end{equation} for all $n\geq 1$. From (\ref{4.12}) and (\ref{4.13}) we
have that
$$ {p_n r_n^{-1} z_n^{-1}\over \mu(r_n) }\leq {7\over n+k_0}  {1\over
\mu(r_n)}\leq  {7\over
r_n(n+k_0)} {r_n\over \mu(r_n)}
 \leq 7\;{\mu({1\over n+k_0})\over  {1\over n+k_0}} {1\over  {\mu(r_n)\over
r_n}}\leq 7.
$$ The proof is complete.

\end{document}